\documentclass[a4paper,12pt]{article}
\usepackage {graphicx}
\usepackage{amsmath}
\usepackage{amssymb}
\usepackage{amsthm}
\usepackage{amscd}
\usepackage{amsfonts}
\usepackage{graphicx}
\usepackage{fancyhdr}
\usepackage[latin1]{inputenc}
\usepackage{tikz}

\newtheorem{thm}{Theorem}

\title {\bf A short proof of Cayley's tree formula}
\author {Alok Shukla}
\date{}
\begin {document}
\maketitle
\begin{center}
\bf {Abstract}
\end{center}
We give a short proof of Cayley's tree formula for counting the number of different labeled trees on $n$ vertices. The following nonlinear recursive relation for the number of labeled trees on $n$ vertices is deduced from a combinatorial argument, 
$$ T_n =  \frac{n}{2}  \sum_{k=0}^{n-2} \left ( \begin {array} {c} n-2 \\ k
\end {array} \right ) T_{k+1} T_{n-k-1}; \ \  for \  n > 1 \ and \  T_1 = 1, $$
and then it is proved that $T_n = n^{n-2}$, which gives yet another proof of the celebrated Cayley's tree formula. 
\section {Introduction}
Cayley's tree formula is a very elegant result in Graph Theory. The problem is to find the number of all possible trees on a given set of labeled vertices. For $ n=2 $ and vertex set $\bf \{v_1,v_2\} $, we have only one tree. For $ n=3 $ and vertex set $\bf \{v_1,v_2,v_3\} $, we have $ 3 $ different trees. Similarly for $ n=4 $, we have $ 16 $ trees. Refer \cite{Gunter} for more details and also for several different proofs of Cayley's tree formula. 
\begin{figure}
\centering
\begin{minipage}{45ex}
\begin{tikzpicture}[scale = 0.6]
\draw [very thick] (3,0)  -- (4,0);
\fill [red] (3,0) circle (0.1);
\fill [blue] (4,0) circle (0.1);
\node at (3.5,0.5) [above] {Labeled trees on two vertices, $ T_2 =1 $.};
\draw [very thick] (0,-2)  -- (1,-2);
\draw [very thick] (0,-2)  -- (.5,-2.87);
\fill [red] (0,-2) circle (0.1);
\fill [blue] (1,-2) circle (0.1);
\fill [green] (.5,-2.87) circle (0.1);

\draw [very thick] (3,-2)  -- (4,-2);
\draw [very thick] (4,-2)  -- (3.5,-2.87);
\fill [red] (3,-2) circle (0.1);
\fill [blue] (4,-2) circle (0.1);
\fill [green] (3.5,-2.87) circle (0.1);

\draw [very thick] (6,-2)  -- (6.5,-2.87);
\draw [very thick] (7,-2)  -- (6.5,-2.87);
\fill [red] (6,-2) circle (0.1);
\fill [blue] (7,-2) circle (0.1);
\fill [green] (6.5,-2.87) circle (0.1);
\node at (3.5,-1.5) [above] {Labeled trees on three vertices, $ T_3 =3 $.};

\draw [very thick] (0,-6)  -- (1,-6);
\draw [very thick] (0,-5)  -- (0,-6);
\draw [very thick] (1,-5)  -- (1,-6);
\fill [red] (0,-5) circle (0.1);
\fill [blue] (1,-5) circle (0.1);
\fill [green] (0,-6) circle (0.1);
\fill [yellow] (1,-6) circle (0.1);
\node at (3.5,-4.5) [above] { $  $ Labeled trees on four vertices, $ T_4 =16 $.};

\draw [very thick] (2,-5)  -- (3,-5);
\draw [very thick] (2,-6)  -- (3,-6);
\draw [very thick] (2,-5)  -- (2,-6);
\fill [red] (2,-5) circle (0.1);
\fill [blue] (3,-5) circle (0.1);
\fill [green] (2,-6) circle (0.1);
\fill [yellow] (3,-6) circle (0.1);

\draw [very thick] (4,-5)  -- (5,-5);
\draw [very thick] (4,-5)  -- (4,-6);
\draw [very thick] (5,-5)  -- (5,-6);
\fill [red] (4,-5) circle (0.1);
\fill [blue] (5,-5) circle (0.1);
\fill [green] (4,-6) circle (0.1);
\fill [yellow] (5,-6) circle (0.1);

\draw [very thick] (6,-5)  -- (7,-5);
\draw [very thick] (6,-6)  -- (7,-6);
\draw [very thick] (7,-5)  -- (7,-6);
\fill [red] (6,-5) circle (0.1);
\fill [blue] (7,-5) circle (0.1);
\fill [green] (6,-6) circle (0.1);
\fill [yellow] (7,-6) circle (0.1);

\draw [very thick,red] (0,-7)  -- (1,-7);
\draw [very thick,] (0,-7)  -- (1,-8);
\draw [very thick,] (1,-7)  -- (0,-8);
\fill [red] (0,-7) circle (0.1);
\fill [blue] (1,-7) circle (0.1);
\fill [green] (0,-8) circle (0.1);
\fill [yellow] (1,-8) circle (0.1);

\draw [very thick,] (3,-7)  -- (3,-8);
\draw [very thick,] (2,-7)  -- (3,-8);
\draw [very thick,] (3,-7)  -- (2,-8);
\fill [red] (2,-7) circle (0.1);
\fill [blue] (3,-7) circle (0.1);
\fill [green] (2,-8) circle (0.1);
\fill [yellow] (3,-8) circle (0.1);

\draw [very thick,] (4,-8)  -- (5,-8);
\draw [very thick,] (4,-7)  -- (5,-8);
\draw [very thick,] (5,-7)  -- (4,-8);
\fill [red] (4,-7) circle (0.1);
\fill [blue] (5,-7) circle (0.1);
\fill [green] (4,-8) circle (0.1);
\fill [yellow] (5,-8) circle (0.1);

\draw [very thick,] (6,-7)  -- (6,-8);
\draw [very thick,] (6,-7)  -- (7,-8);
\draw [very thick,] (7,-7)  -- (6,-8);
\fill [red] (6,-7) circle (0.1);
\fill [blue] (7,-7) circle (0.1);
\fill [green] (6,-8) circle (0.1);
\fill [yellow] (7,-8) circle (0.1);

\draw [very thick,] (0,-9)  -- (1,-9);

\draw [very thick,] (1,-9)  -- (0,-10);
\draw [very thick,] (0,-10)  -- (1,-10);
\fill [red] (0,-9) circle (0.1);
\fill [blue] (1,-9) circle (0.1);
\fill [green] (0,-10) circle (0.1);
\fill [yellow] (1,-10) circle (0.1);

\draw [very thick,] (2,-9)  -- (3,-9);

\draw [very thick,] (3,-10)  -- (2,-9);
\draw [very thick,] (2,-10)  -- (3,-10);
\fill [red] (2,-9) circle (0.1);
\fill [blue] (3,-9) circle (0.1);
\fill [green] (2,-10) circle (0.1);
\fill [yellow] (3,-10) circle (0.1);

\draw [very thick,] (4,-9)  -- (4,-10);

\draw [very thick,] (4,-10)  -- (5,-9);
\draw [very thick,] (5,-9)  -- (5,-10);
\fill [red] (4,-9) circle (0.1);
\fill [blue] (5,-9) circle (0.1);
\fill [green] (4,-10) circle (0.1);
\fill [yellow] (5,-10) circle (0.1);

\draw [very thick,] (6,-9)  -- (6,-10);
\draw [very thick,] (6,-9)  -- (7,-10);
\draw [very thick,] (7,-9)  -- (7,-10);

\fill [red] (6,-9) circle (0.1);
\fill [blue] (7,-9) circle (0.1);
\fill [green] (6,-10) circle (0.1);
\fill [yellow] (7,-10) circle (0.1);


\draw [very thick,] (0,-11)  -- (1,-11);

\draw [very thick,] (0,-11)  -- (0,-12);
\draw [very thick,] (0,-11)  -- (1,-12);

\fill [red] (0,-11) circle (0.1);
\fill [blue] (1,-11) circle (0.1);
\fill [green] (0,-12) circle (0.1);
\fill [yellow] (1,-12) circle (0.1);

\draw [very thick,] (2,-11)  -- (3,-11);

\draw [very thick,] (3,-11)  -- (3,-12);
\draw [very thick,] (2,-12)  -- (3,-11);

\fill [red] (2,-11) circle (0.1);
\fill [blue] (3,-11) circle (0.1);
\fill [green] (2,-12) circle (0.1);
\fill [yellow] (3,-12) circle (0.1);

\draw [very thick,] (5,-11)  -- (5,-12);

\draw [very thick,] (4,-11)  -- (5,-12);
\draw [very thick,] (4,-12)  -- (5,-12);

\fill [red] (4,-11) circle (0.1);
\fill [blue] (5,-11) circle (0.1);
\fill [green] (4,-12) circle (0.1);
\fill [yellow] (5,-12) circle (0.1);

\draw [very thick,] (6,-11)  -- (6,-12);

\draw [very thick,] (6,-12)  -- (7,-11);
\draw [very thick,] (6,-12)  -- (7,-12);
\fill [red] (6,-11) circle (0.1);
\fill [blue] (7,-11) circle (0.1);
\fill [green] (6,-12) circle (0.1);
\fill [yellow] (7,-12) circle (0.1);
\node at (7,-13) {};
\end{tikzpicture}
\end{minipage}
\caption{Cayley's tree.}
\end{figure}
\section {Counting trees}
\begin{thm}
\emph{There are exactly $\bf n^{n-2}$ labeled trees on
$ n $ vertices.}
\end{thm}
\noindent \textbf {Proof:}    
Let $\bf T$ represent the set of all labeled trees on $\bf n$ vertices of the given vertex set  $\bf V\{v_1,v_2 \ldots v_n \}$. Let $\bf T_n$ denote the cardinality of $\bf T$. Let $\bf E_n$ denote the number of trees in $\bf T$ that always contain one specified edge, say $\bf v_1v_2$. It should be noted that $\bf T $ is not affected by any rearrangement of vertices of $\bf V$. This symmetry implies that all the edges are equivalent and contribute equally to $\bf T $. As the number of possible edges is $ {{n} \choose {2}} = \frac{1}{2}n(n-1)$; the total number of edges in all trees in $\bf T$ equals
\begin {equation*}
\bf \frac {1} {2} \ n \ (n - 1) E_n = (n-1) T_n \,,
\end {equation*}
where the right side of equality follows from the fact that each one of $\bf T_n$ trees in $\bf T$ contains $(n-1)$ edges.
Therefore we have for all $\bf n>1$,
\begin {equation}
\label {Eq1}
\bf T_n = \frac {1} {2} \ n \ E_n \,.
\end{equation}
\noindent Further, it is easy to see that $\bf T_1= 1 $. 
Without any loss of generality, we fix edge $\bf v_1v_2 $ and count the number of possible trees to determine $\bf E_n$. For this, $\bf k$ vertices are chosen from the set $\bf V\{v_3,v_4,\ldots v_n \}$ in $\bf n-2 \choose k $ ways. These $\bf k $ vertices together with vertex $\bf v_1$ form a set of $\bf (k+1)$ vertices and would give $\bf T_{k+1}$ trees. The remaining $\bf (n-2-k) $ vertices are attached to vertex $\bf v_2$ to from a set of $\bf (n-k-1)$ vertices to yield $\bf T_{n-k-1} $ trees. By summing for $\bf k=0 $ to $\bf n-2$ we get,
\begin {equation}
\label {Eq2}
\bf  E_n = \sum _{\bf k=0}^{\bf n-2} {{n-2} \choose {k}} T_{k+1} T_{n-k-1} \,.
\end {equation}
From Eq(\ref{Eq1}) and Eq(\ref{Eq2}) we get,
\begin {equation}
\label {Eq3}
\bf T_n =  \frac{1}{2} n \sum _{\bf k=0}^{\bf n-2} {{n-2} \choose {k}} T_{k+1} T_{n-k-1} \,.
\end {equation}
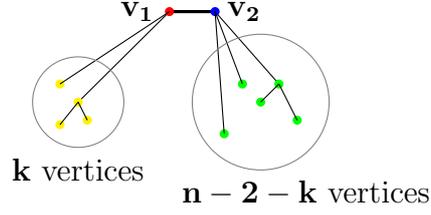
\begin{figure}
\centering
\begin{minipage}{45ex}
\begin{tikzpicture}[scale = 0.6]
\draw [very thick] (3,0)  -- (4,0) node [right] {$ \bf v_2 $};
\node at (2.9,0) [left] {$\bf v_1 $};
\draw [-] (1,-2)--(3,0) ;
\draw [-] (5.4,-1.6)--(4,0) ;
\fill [red] (3,0) circle (0.1);
\fill [blue] (4,0) circle (0.1);
\fill [yellow] (1,-2) circle (0.1);
\fill [green] (5,-2) circle (0.1);
\fill [yellow] (1.2,-2.4) circle (0.1);
\fill [green] (5.8,-2.4) circle (0.1);
\fill [green] (5.4,-1.6) circle (0.1);
\fill [yellow] (0.6,-1.6) circle (0.1);
\fill [yellow] (0.6,-2.5) circle (0.1);
\fill [green] (4.6,-1.6) circle (0.1);
\fill [green] (4.2,-2.7) circle (0.1);
\draw [-] (4.2,-2.7)--(4,0) ;
\draw [-] (4.6,-1.6)--(4,0) ;
\draw [-] (5.4,-1.6)--(5.8,-2.4) ;
\draw [-] (5.4,-1.6)--(5,-2) ;
\draw [-] (1.2,-2.4) -- (1,-2);
\draw [-] (0.6,-1.6) -- (3,0);
\draw [-] (0.6,-2.5) -- (1,-2);
\draw [gray] (1,-2) circle(1);
\draw [gray] (5,-2) circle(1.5);
\node at (1,-3) [below] {$\bf k $ vertices};
\node at (6,-3.5) [below] {$\bf n-2- k $ vertices};
\end{tikzpicture}
\end{minipage}
\caption{Counting all the tree that contain the edge $ \bf v_1v_2 $.}
\end{figure}
\noindent Let the exponential generating function $\bf T(S)$ be defined as
\begin {align}
\label {Eq4}
\ \bf T(S) = \sum_{n=1}^{\infty} T_n \frac{S^n}{(n-1)!} \,.
\end {align}
Then,
\begin {align*}
\bf (T(S))^{2} &=  \bf \sum_{n=1}^{\infty} \sum_{k=0}^{n-2} \left \{ T_{k+1} \frac {S^{k+1}}{(k!)} \right \} \left \{T_{n-k-1} \frac {S^{n-k-1}} {(n-k-2)!}\right \} \\ 
\bf &= \bf \sum_{n=1}^{\infty} \sum_{k=0}^{n-2} \left \{ {{n-2} \choose {k}}  T_{k+1} T_{n-k-1} \frac {S^{n}} {(n-2)!}\right \}.
\end {align*}
\noindent By using the value of $\bf T_n$ from the recurrence relation of Eq (\ref{Eq3}), 
\begin {equation} 
\label {Eq5}
\bf (T(S))^2 = \sum_{n=1}^{\infty} \left \{ 2 \frac{(n-1)}{n} T_n \frac {S^{n}} {(n-1)!}\right  \}. 
\end {equation}
\noindent On differentiating with respect to $\bf S$, 
\begin {equation*}
\bf  2 T(S)T'(S)  = \sum_{n=1}^{\infty} \left \{ 2 (n-1) T_n \frac {S^{n-1}} {(n-1)!}\right \}.
\end{equation*} 
Therefore,
\begin {equation*}
\bf T(S)T'(S) = \sum_{n=1}^{\infty} \left \{ n T_n \frac {S^{n-1}} {(n-1)!}\right \} - 
\frac{1}{S} \sum_{n=1}^{\infty} \left \{T_n \frac {S^{n}} {(n-1)!}\right \} = T'(S) - \frac {T(S)} {S}.
\end {equation*}
So we get,
\begin {equation*}
\bf \partial [T(S)] = \partial[ln \left [\frac {T(S)}{S}\right ]  ].
\end {equation*}
\noindent On integrating,
\begin {equation}
\label {Eq6}
\bf T(S) = ln \left [\frac {T(S)}{S}\right ] + C \, ,
\end {equation}
where $ C $ is a constant of integration.
\\
To fix $ C $, put $ S=0 $ in Eq(\ref{Eq6}) and 
note that from Eq(\ref{Eq4}), $T(0) = 0$ and $\frac {T(S)}{S}|_{S=0} = T_1 = 1$. So we have 
$ 0 = ln(1) + C \Rightarrow C = 0$ 
and Eq(\ref{Eq6}) reduces to 
\begin {equation}
\label {Eq7}
\bf T(S) = S \exp^{\bf T(S)}.
\end {equation}
\noindent One can use \textit{\textbf{Lagrange Inversion Theorem}} (see \cite{Flajolet}) to expand $T(S)$ in power series of $S$
and the result is 
\begin {equation}
\label {Eq8}
\bf T(S)  =  \bf \sum_{n=1}^{\infty} \frac {n^{n-1}}{n!}S^n = \bf \sum_{n=1}^{\infty} n^{n-2}\frac {S^n}{(n-1)!} \,.
\end {equation}
\noindent By equating the coefficient of $S^n$ in Eq (\ref {Eq4}) and Eq (\ref {Eq8}) we get $\bf T_n  =  \bf n^{n-2}$. 
\newline This completes the proof.
\section {Conclusion}
A number of remarkable proofs of Cayley's tree formula are known. The proof presented in this paper exploits the symmetry of edges. It is quite interesting to witness different areas of mathematics converging so beautifully to yield yet another proof of this celebrated result.

\end {document}